\documentclass[11pt,tbtags,reqno]{article}
\usepackage[margin=0.96in]{geometry}
\usepackage{amssymb}
\usepackage{amsmath}
\usepackage{amsthm}
\usepackage[swedish, english]{babel}
\usepackage{psfrag}   
\usepackage{graphicx}
\usepackage[retainorgcmds]{IEEEtrantools}
\usepackage{bbm}
\usepackage[square, comma, numbers, sort&compress]{natbib}
\usepackage{enumerate}
\usepackage{soul}
\usepackage{hyperref}
\usepackage[colorinlistoftodos, textsize=tiny, backgroundcolor=yellow
, linecolor=green]{todonotes}
\usepackage{sidecap}
\usepackage{float}
\usepackage{leftidx}

\usepackage{bm}

\newcommand{\E}{\mathbb E}
\newcommand{\R}{\mathbb R}
\newcommand{\N}{\mathbb N}

\numberwithin{equation}{section}
\theoremstyle{plain}

\newtheorem{theorem}{Theorem}[section]
\newtheorem{lemma}[theorem]{Lemma}

\newtheorem{remark}[theorem]{Remark}
 \newtheorem{example}[theorem]{Example}
\theoremstyle{definition}
 \theoremstyle{example}

\begin{document}

\title{\textsc{
A   Weak Law of Large Numbers for Dependent Random Variables}  \thanks{~ We are indebted to Albert \textsc{Shiryaev} for pointing out to us the concept of  generalized mathematical expectations  in \textsc{Kolmogorov} \cite{Kol}, which prompted this investigation.  We are  grateful to Nathan \textsc{Soedjak} for his careful reading of the manuscript and his incisive comments. 
}
}

\author{  
\textsc{Ioannis Karatzas} \thanks{~
Departments of Mathematics and Statistics,  Columbia University, New York, NY 10027 (e-mail: {\it ik1@columbia.edu}). Support from   National Science Foundation   Grant  DMS-20-04977 is gratefully acknowledged.  
}  
 \and
\textsc{Walter Schachermayer}                \thanks{  ~
Faculty of Mathematics, University of Vienna, Oskar-Morgenstern-Platz 1, 1090 Vienna, Austria (email: {\it walter.schachermayer@univie.ac.at}). Support from the  Austrian Science Fund (FWF) under
 grant P-28861 and grant P-35197 is gratefully acknowledged.    
          }
                                      }

\maketitle
\begin{abstract}
\noindent
Every sequence $f_1, f_2, \cdots \, $ of random variables with $\, \lim_{M \to \infty } \big( M \, \sup_{k \in \N} \mathbb{P} ( |f_k| > M ) \big)=0\,$ contains a subsequence $  f_{k_1}, f_{k_2} , \cdots  \,$ that satisfies, along with all its subsequences, the weak law of large numbers: $ \, \lim_{N \to \infty} \big( (1/N) \sum_{n=1}^N f_{k_n} - D_N \big) =0\,,$  in probability. Here   $\,  D_N\,   $  is a   ``corrector" random variable  with  values in  $[-N,N]$, for each $\,   N \in \N \,$; these correctors are all equal to zero if, in addition, $\, \liminf_{k \to \infty}  \mathbb{E} \big( f_k^2 \, \mathbf{ 1}_{ \{ |f_k| \le M \} } \big) =0\,$ holds for every $M \in (0, \infty)\,.$ 
\end{abstract}



\bigskip


\section{Introduction}
\label{sec1}


 On a probability space $(\Omega, \mathcal{F}, \mathbb{P})$,   consider    real-valued   measurable functions $f_1, f_2, \cdots \,.$ If these are independent and have  the same distribution with     $\E ( | f_1|) < \infty\,$, the celebrated \textsc{Kolmogorov}  strong law of large numbers (SLLN: \cite{Kol};\,\cite{K2};\,\cite{Du},\,section\,2.4) states that the ``sample average" $ \, ( f_1 + \cdots + f_N) / N\,$ converges $\mathbb{P}-$a.e.\,to the ``ensemble average" $\,\E (f_1) = \int_\Omega f_1 \, \mathrm{d} \mathbb{P}\,,$ as $ N \to \infty$. 
 
  \smallskip
 A bit more generally,  if the functions $f_k (\omega) = f \big( T^{k-1}(\omega) \big), \, k \ge 2, \, \omega \in \Omega$ are the images of an integrable function $f_1: \Omega \to \R$ along the orbit  of successive actions of a measure-preserving transformation $T:   \Omega \to \Omega\,,$ then the  above sample average converges $\mathbb{P}-$a.e.\,as $N \to \infty$ to the conditional expectation  $f_* = \E ( f_1 |{\cal I})\,$  of $f_1$, given    the $\sigma-$algebra  $\,{\cal I}\,$   of   $\,T-$invariant sets,  by   \textsc{Birkhoff}'s pointwise ergodic theorem (\cite{B};\,\cite{Du}, p.\,333).

 \smallskip
A deep  result of \textsc{Koml\'os} \cite{Kom}, already   55 years  old but always very striking, says that such ``stabilization via  averaging" occurs within {\it any} sequence $f_1, f_2, \cdots \,$ of measurable, real-valued  functions with  $\, \sup_{n \in \N} \E ( | f_n|) < \infty\,.$ More precisely, there exist  then an integrable function $f_*$ and a subsequence $\big\{ f_{k_n} \big\}_{n \in \N}\,$ such that   $ \, ( f_{k_1} + \cdots + f_{k_N}) / N\,$ converges to $f_*\,$, $\,\mathbb{P}-$a.e.\,\,as $ N \to \infty$; and   the same is true   ``hereditarily", that is, for any further subsequence  of $\big\{ f_{k_n} \big\}_{n \in \N}\,.$  

\smallskip
 We have also another celebrated result of \textsc{Kolmogorov},  the weak  law of large numbers (WLLN: \cite{Kol};\,\cite{Ch},\,section 5.2;\,\cite{Du},\,\S\,2.2.3) for a sequence $f_1, f_2, \cdots \,$ of real-valued,   measurable functions which are independent. If these  are identically  distributed and satisfy the weak-$\mathbb{L}^1-$type condition 
\begin{equation}
\label{1.1}
\lim_{M \to \infty} \, \Big( M \cdot \mathbb{P} \big( | f_1 | > M \big) \Big) \,=\,0
\end{equation} 
(rather than $\E ( | f_1|) < \infty\,$), then we have the WLLN
\begin{equation}
\label{1.2}
\lim_{N \to \infty} \left( \frac{1}{\,N\,} \sum_{n=1}^N f_n - D_N \right)  \,=\,0\,, \quad \text{in probability}
\end{equation} 
for the sequence of ``correctors"  
\begin{equation}
\label{1.3}
D_N\, :=\, \mathbb{E} \Big( f_1\, \mathbf{ 1}_{\{ | f_1| \le N  \} } \Big)\,, \quad N \in \N\,;
\end{equation} 
whereas, if the independent functions $f_1, f_2, \cdots \,$ do not have the same distribution but satisfy 
\begin{equation}
\label{1.4}
\lim_{N \to \infty} \sum_{n=1}^N \, \mathbb{P} \big( | f_n | > N \big) \,=\,0\,, \qquad \lim_{N \to \infty} \frac{1}{N^2} \sum_{n=1}^N \, \mathbb{E} \Big( f_n^2\, \,\mathbf{ 1}_{\{ | f_n| \le N  \} } \Big)\,=\,0\,,
\end{equation}
then again the convergence in probability (WLLN) in ({\ref{1.2}) holds, though now with   correctors
\begin{equation}
\label{1.5}
D_N\, :=\, \frac{1}{N } \sum_{n=1}^N \, \mathbb{E} \Big( f_n\, \mathbf{ 1}_{\{ | f_n| \le N  \} } \Big)\,, \quad N \in \N\,.
\end{equation}
It was shown in  \cite{GK},\,\cite{F},\,(\cite{Ch},\,Theorem 5.2.3)  that, for independent   $f_1, f_2, \cdots \,,$  the conditions in (\ref{1.4}) are not only sufficient but also necessary for the existence of a sequence $\,D_1, D_2, \cdots \,$ of real numbers with the property (\ref{1.2}). 

Let us also note, that   the correctors in both (1.3), (1.5) satisfy $\, | D_N| \le N\,;$ and that they are all equal to zero, if each of the $f_1, f_2, \cdots \,$  has distribution   symmetric around the origin. 

\smallskip
The purpose of this Note is to present a \textsc{Koml\'os}-type version of the weak law of large numbers. This is formulated in the next section, and proved in section \ref{sec3}. The proof, considerably simpler than its counterpart for the strong law in \cite{Kom}, is based on truncation and on weak-$\mathbb{L}^2$ convergence arguments, which give also sufficient conditions for the resulting correctors to be equal to zero. Examples and ramifications   are taken up in section \ref{sec4}.


\section{Result}
\label{sec2}


   We consider    real-valued   measurable functions $f_1, f_2, \cdots \,$  on a probability space $(\Omega, \mathcal{F}, \mathbb{P})$, and introduce for every $M \in (0, \infty) $  the quantities
   \begin{equation}
\label{2.1}
\tau_n (M)\, := \,M \cdot  \mathbb{P} \big( | f_n | > M \big)\,, \qquad \tau (M) \,:=\, \sup_{n \in \N} \, \tau_n (M)\,. 
\end{equation}

\begin{theorem} 
\label{thm2.1}
{\bf A General,  Hereditary WLLN.} In the above context, suppose that the weak-$\mathbb{L}^1-$type condition 
   \begin{equation}
\label{2.2}
\lim_{M \to \infty} \tau (M) \,=\, 0  
\end{equation}
holds. There exist then a sequence of ``corrector" random variables $D_1, D_2, \cdots $ with 
  \begin{equation}
\label{2.3}
 \mathbb{P} \big( | D_N | \le N \big) \,=\, 1 \qquad \text{for every} ~~ N \in \N\,,
\end{equation}
and 
a subsequence $\big\{ f_{k_n} \big\}_{n \in \N}\,, $ such that the WLLN   
 \begin{equation}
 \label{2.4}
 \, \lim_{N \to \infty} \Big( \frac{1}{N} \sum_{n=1}^N f_{k_n} - D_N \Big) =0\,, \quad \hbox{ in probability} 
 \end{equation}
 is satisfied    ``hereditarily"; i.e., not just by $\big\{ f_{k_n} \big\}_{n \in \N}\,  $ but also by  all its subsequences. 
If, in addition, 
\begin{equation}
\label{2.5}
 \liminf_{k \to \infty} \, \mathbb{E} \Big( f_k^2 \,\, \mathbf{ 1}_{ \{ |f_k| \le M \} } \Big) \,=\,0
 \end{equation}
  holds for every $M \in (0, \infty)\,,$ the correctors in  \eqref{2.4}, \eqref{2.3}  can be chosen as   $ D_N =0 $ for every $N \in \N.$ 
\end{theorem} 

The correctors $D_1, D_2, \cdots $  
correspond to the   generalized mathematical expectations in \textsc{Kolmogorov} \cite{Kol},\,\S 6.4; they are also related to the  nonlinear expectations  developed by \textsc{Peng} in \cite{P}.


\section{Proof}
\label{sec3}

We start with the simple but crucial idea of ``truncation". This  goes back to the work of  \textsc{Khintchine} and \textsc{Kolmogorov} (\cite{GK},\,\cite{K1}), where it plays a major role in the proofs of   laws of large numbers and of convergence results for series of random variables. 

\begin{lemma}
Under the condition \eqref{2.2}, we have 
\begin{equation}
\label{3.1}
\, \lim_{N \to \infty} \bigg( \frac{1}{N} \sum_{n=1}^N f_{n}
\,- \, \frac{1}{N} \sum_{n=1}^N f_{n}\,\, \mathbf{ 1}_{ \{ |f_n| \le N \} }
\bigg) =0\,, \quad \hbox{ in probability}.
\end{equation}
\end{lemma}

\noindent
{\it Proof:} For every $\varepsilon >0\,,$ the expression 
$$
 \mathbb{P} \left( \bigg| \frac{1}{N} \sum_{n=1}^N f_{n}\,\, \mathbf{ 1}_{ \{ |f_n| > N \} } \bigg| > \varepsilon \right) \,\le \,  \mathbb{P} \left( \bigcup_{n=1}^N \big\{ |f_n| > N \big\} \right)  \,\le \, \sum_{n=1}^N  \mathbb{P} \big( |f_n| > N \big) \,\le \, N \cdot \max_{1 \le n \le N} \mathbb{P} \big( |f_n| > N \big) 
 $$ 
is dominated by $\, N \cdot \sup_{n \in \N} \mathbb{P} \big( |f_n| > N \big)$, which converges to zero as $ N \uparrow \infty$ on account of   (\ref{2.2}). \qed

\medskip
It follows that, in order to establish (\ref{2.4}), it is enough to prove 
\begin{equation}
\label{3.2}
\, \lim_{N \to \infty} \bigg( \frac{1}{N} \sum_{n=1}^N f_{n}\, \mathbf{ 1}_{ \{ |f_n| \le N \} } - D_N \bigg) =0\,, \quad \hbox{ in probability} 
\end{equation}
for a suitable sequence $D_1, D_2, \cdots $ of correctors, and 
along an appropriate subsequence of $\{ f_n \}_{n \in \N}\,$ denoted by the same symbols for economy of exposition---as well as along all further subsequences of this subsequence.  

\medskip
\noindent
\textsc{Proof of Theorem \ref{thm2.1}:} For each integer $N \in \N$ we consider the truncated functions
\begin{equation}
\label{3.3}
f^{\,[-N,N]}_n \,:=\, f_{n}\, \mathbf{ 1}_{ \{ |f_n| \le N \} }\,, \qquad n \in \N
\end{equation}
that appear in (\ref{3.1}), (\ref{3.2}). These are bounded in $\mathbb{L}^\infty$ (as they take values in  $[-N, N]$), thus   bounded in $\mathbb{L}^2$  as well. As a result we can extract,  for each $ N \in \N\,,$  a subsequence of $\{ f_n \}_{n \in \N}\,$ denoted by the same symbols for economy of exposition,   such that  the sequence in (\ref{3.3}) converges weakly in $\mathbb{L}^2$ to some $ D_N \in  \mathbb{L}^2\,$: 
\begin{equation}
\label{3.4}
\lim_{n \to \infty}\, 
\mathbb{E} \Big( f^{\,[-N,N]}_n \cdot \xi \Big) \,=\,\mathbb{E} \big( D_N \cdot \xi \big)\,, \qquad \forall~~ \xi  \in \mathbb{L}^2\,.
\end{equation}
And by standard 
diagonalization arguments, we can extract then  a further subsequence of $\{ f_n \}_{n \in \N}\,,$ denoted again by the same symbols, and such that \eqref{3.4} holds for {\it every} $N \in \N$. 

It is fairly straightforward to check that these weak-$\mathbb{L}^2$ limits in \eqref{3.4} satisfy (\ref{2.3}).  
On the other hand,  the lower-semicontinuity of the $\mathbb{L}^2-$norm under weak-$\mathbb{L}^2$ convergence, in this case 
$$
\big \| D_N \big \|_{\mathbb{L}^2} \, \le \, \liminf_{n \to \infty} \Big \| f^{\,[-N,N]}_n \Big \|_{\mathbb{L}^2}\,,
$$
 gives $\mathbb{P} (D_N=0)=1$ for every $N \in \N$, under the condition (\ref{2.5}).  

\smallskip
We introduce now, for each $M \in (0, \infty),$  the quantities 
\begin{equation}
\label{3.5}
\sigma_n (M) \,:=\, \frac{1}{M} \, \,\mathbb{E} \Big( f_n^2 \,\, \mathbf{ 1}_{ \{ |f_n| \le M \} } \Big)\,,  \qquad \sigma  (M) \,:=\, \sup_{n \in \N} \, \sigma_n (M) \,. 
\end{equation}
As shown by \textsc{Feller} (\cite{F},\,p.\,235; see also \cite{Du},\,\S\,2.3.3), these quantities are related to those in (\ref{2.1}) via
\begin{equation}
\label{3.6}
0 \le \sigma_n (M)\,=\, \frac{2}{M} \int_0^M \tau_n (t) \, \mathrm{d} t - \tau_n (M)  \,\le\,  \frac{2}{M} \int_0^M \tau  (t)\, \mathrm{d} t
\end{equation} 
for every $n \in \N, \, M \in (0, \infty)$, \footnote{~   In the integrand of this expression, as it appears on page 235 of \cite{F}, there is a typo; this is here corrected. } thus 
\begin{equation}
\label{3.7}
0 \le \sigma  (M)\,\le\, \frac{2}{M} \int_0^M \tau  (t) \, \mathrm{d} t \,, \qquad M \in (0, \infty)\,.
\end{equation}

From this bound (\ref{3.7}) and the assumption (\ref{2.2}), it follows that we have
\begin{equation}
\label{3.8}
\lim_{M \to \infty} \sigma (M) \,=\, 0\,. 
\end{equation}
We note also 
\begin{equation}
\label{3.9}
\mathbb{E}  \Big( f^{\,[-M,M]}_n     \Big)^2 \,=\,
\mathbb{E} \Big( f_n^2 \,\, \mathbf{ 1}_{ \{ |f_n| \le M \} } \Big) \,\le\, M \cdot  \sigma_n (M) \,\le\,M \cdot  \sigma (M)
\end{equation}
for all $n \in \N\,,$ $ M \in (0, \infty)\,,$ therefore
\begin{equation}
\label{3.10}
\mathbb{E} \big( D_M^2 \big)  \,\le\,
\sup_{n \in \N} \,\mathbb{E}  \Big( f^{\,[-M,M]}_n     \Big)^2   \,\le\,M \cdot  \sigma (M) \,=\, o (M)\,, \quad \text{as}~ M \to \infty\,.
\end{equation}

We observe at this point that, in order to prove (\ref{3.2}) and thus (\ref{2.4}) as well, along a suitable subsequence, it is enough to show convergence along this subsequence in $\mathbb{L}^2$, namely
\begin{equation}
\label{3.11}
\lim_{N \to \infty} \frac{1}{N^2} \cdot 
\mathbb{E} \left( \sum_{n=1}^N \Big( f^{\,[-N,N]}_n - D_N   \Big)
\right)^2 \,=\,0\,.
\end{equation}
And developing the square, we need to show that the expectations of both the sum of squares
\begin{equation}
\label{3.12}
 \sum_{n=1}^N \mathbb{E}  \Big( f^{\,[-N,N]}_n - D_N   \Big)^2 \,\le\, 
2 \,\sum_{n=1}^N \,\mathbb{E}  \Big( f^{\,[-N,N]}_n   \Big)^2 + 2\,N \cdot \mathbb{E} \big( D_N \big)^2
\end{equation} 
and the sum  of   cross-products 
\begin{equation}
\label{3.13}
2 \, \sum_{n=1}^N \sum_{1 \le j <n}  
 \mathbb{E} \left[  \Big( f^{\,[-N,N]}_j - D_N   \Big) \Big( f^{\,[-N,N]}_n - D_N   \Big) \right] 
\end{equation}
are of order $o (N^2)$, as $N \to \infty\,,$ for the subsequence in question and for all its subsequences. Now,  from (\ref{3.9}), (\ref{3.10}), the expression in (\ref{3.12}) is already dominated by $\, 4\,N^2 \cdot \sigma (N) = o (N^2),$ as $N \to \infty$, on account of (\ref{3.8}). 

\smallskip
Let us recall what happens at this juncture in the case of independent $f_1, f_2, \cdots \,:$ the correctors $D_N$ are real constants, given as in (\ref{1.5}), so the differences  $f^{\,[-N,N]}_n - D_N\,, ~ n=1, \cdots ,N\,$ are independent with zero mean, thus uncorrelated. The expectations of their cross-product in (\ref{3.13}) vanish, and the argument ends here. 

In the general case, when {\it nothing} is assumed about the finite-dimensional distributions of the $f_1, f_2, \cdots $ (in particular, when these functions are not independent)  we need to guarantee, by passing to a further subsequence if necessary,  that the expression in (\ref{3.13}) is also of order $o (N^2)$, as $N \to \infty\,.$ One way to accomplish  this, is to ensure that the differences  $f^{\,[-N,N]}_n - D_N\,, ~ n=1, \cdots, N\,$  are {\it  very close to being uncorrelated.}

 \smallskip
We do this by induction, in the following  manner: Suppose the terms $\,f_1, \cdots, f_{n-1}\,$ of the subsequence have been chosen. We select the next term $f_n$ in such a way, that the difference $f^{\,[-N,N]}_n - D_N\,,$ with $\, N \le e^{\,n^2},$ is ``almost orthogonal" to all of the $f^{\,[-N,N]}_1 - D_N\,,\cdots, f^{\,[-N,N]}_{n-1} - D_N\,;$ to wit, 
\begin{equation}
\label{3.14}
\bigg| \, \mathbb{E} \left[  \Big( f^{\,[-N,N]}_j - D_N   \Big) \Big( f^{\,[-N,N]}_n - D_N   \Big) \right] \bigg| \,\le\, e^{-n^2} \, \le \, \frac{1}{N}
\end{equation}
for every $\, j = 1, \cdots, n-1$. Such a choice of $f_n$ is certainly possible on account of (\ref{3.4}), and completes the induction step. 

Returning to (\ref{3.13}), we note that the summation
$$
2 \, \sum_{n=1}^{\sqrt{ \,\log N\,} } \sum_{1 \le j <n} \, \bigg| \,
 \mathbb{E} \left[  \Big( f^{\,[-N,N]}_j - D_N   \Big) \Big( f^{\,[-N,N]}_n - D_N   \Big) \right] \bigg|
$$
is straightforward to control: each summand is bounded 
by $\, N \cdot \sigma (N)\,$ on account of (\ref{3.9}), (\ref{3.10}), so the entire summation is   
of the order 
$$\, 
N  \sigma (N)  \,\sum_{n=1}^{\sqrt{ \,\log N\,} } 2\, n \,\sim \,N  \sigma (N)  \cdot  \log N = o(N^2)\,,
$$ 
as   $N \to \infty$. 
On the other hand, the validity of (\ref{3.14}) for every $\, j = 1, \cdots, n-1$ implies that the summation 
$$
2 \, \sum^{N}_{n=1 + \sqrt{ \,\log N\,} } \,\sum_{1 \le j <n} \, \bigg| \,
 \mathbb{E} \left[  \Big( f^{\,[-N,N]}_j - D_N   \Big) \Big( f^{\,[-N,N]}_n - D_N   \Big) \right] \bigg|
$$
is  of the order 
$$
2 \, \sum^{N}_{n=1 + \sqrt{ \,\log N\,} } \, n\, e^{-n^2/2}\, \sim \, \int_{\sqrt{ \,\log N\,}}^N2\, x\, e^{-x^2/2}\, \mathrm{d} x\,=\, \frac{1}{N} - e^{-N}\,,
$$
as $N \to \infty$, thus certainly of order $\, o (N^2)\,.$ And it follows that the expression of (\ref{3.13}) is of order $\, o (N^2)\,$ as well. 

The argument is now complete. 
It is also straightforward  to check that it works just as well for an arbitrary subsequence, of the subsequence just constructed. \qed



\section{Ramifications and Examples}
\label{sec4}

The condition (\ref{3.2}), which reads
$$\,
\lim_{M \to \infty} \left( \, \sup_{n \in \N} \, \tau_n (M) \right) \,=\,0\,,
 $$
can be weakened to 
\begin{equation}
\label{3.15}
 \lim_{M \to \infty} \left(\, \liminf_{n \in \N} \, \tau_n (M) \right) \,=\,0\,
\end{equation}
Indeed, by passing to a subsequence, this becomes
\begin{equation}
\label{3.16}
 \lim_{M \to \infty} \left(\, \limsup_{n \in \N} \, \tau_n (M) \right) \,=\,0\,,
\end{equation}
and one checks relatively easily that (\ref{3.16}) can replace (\ref{3.2}) in the inductive construction of the subsequence (of) $\big\{ f_n \big\}_{n \in \N}\,.$ {\it We note also that the condition \eqref{3.16} can be satisfied, while \eqref{3.2} fails.} 

\smallskip
To see this, take $\, g \in \mathbb{L}^0$ with 
\begin{equation}
\label{3.17}
 \limsup_{M \to \infty} \Big(\, M \cdot \mathbb{P} \big( |g| > M \big)  \Big) \,>\,0
\end{equation}
and define the functions 
\begin{equation}
\label{3.18}
 f_n \,:=\, g \cdot \mathbf{ 1}_{ \{ |g| >n \} }\,, \quad n \in \N\,.
\end{equation}
We have then $\, \tau_n (M) = M \cdot \mathbb{P} \big( |g| > M \vee n \big) \,,~ \tau (M) = M \cdot \mathbb{P} \big( |g| > M   \big) \,,$ so (\ref{3.17}) means that (\ref{3.2}) fails. However, 
$$
\lim_{n \to \infty }  \tau_n (M)\,=\, M \cdot \lim_{n \to \infty }  \mathbb{P} \big( |g| >   n \big)\,=\,0
$$
holds for every $M \in (0, \infty)$, so (\ref{3.16}) is satisfied. We obtain this way the WLLN (\ref{2.4}) for a suitable sequence of correctors $\,D_1, D_2, \cdots\,.$ 

It is also checked readily that the condition (\ref{2.5}) is   satisfied here, so all these correctors can actually be chosen    equal to  zero. 

\begin{example}   
{\rm
To provide another illustration of Theorem \ref{thm2.1} that highlights the role of both conditions   (\ref{2.1}) and (\ref{2.5}) in a more substantial way, let us revisit  an old example from \cite{K1} (see also section 5.2 of \cite{Ch}) in slightly   modified form.   Suppose that the functions $f_1, f_2, \cdots \,$ satisfy
\begin{equation}
\label{4.1}
 \mathbb{P} \big( f_n =0 \big)    \,=\, \varrho_n\,; \qquad  \mathbb{P} \big( f_n =k \big)    \,=\, \frac{\,(1-\varrho_n) \,c\,}{k^2 \, \log k}\,,~~~ k=2,3, \cdots 
\end{equation}
with constants $0< \varrho_n <1$ and $\, 2\, c = \Big( \sum_{k \ge 2} k^{-2}\, \big(1 /  \log k \big)  \Big)^{-1}\,,$ for every $n \in \N$.

\smallskip
 We do {\it not} impose any condition on the finite-dimensional joint distributions of the $f_1, f_2, \cdots \,$; in particular, we do {\it not} require  the $f_1, f_2, \cdots \,$ to be independent. 

\smallskip
In this setting,   
$$
\tau_n (M) \,= \,2 \,c \,M \big(1-\varrho_n \big) \sum_{k>M} \frac{1}{\, k^2 \log k \,}\,\sim\, \frac{\,2\, c\, }{\,\log M\,} \, \big(1-\varrho_n \big)
$$
holds for integers $M \ge 2$ in the notation of (\ref{2.1}). Thus, 
$\,
\tau (M) \,=\, \sup_{n \in \N} \, \tau_n (M)\,\le\, (  2 \, c )  /  \log M  \,,
$ 
and the condition  (\ref{2.2}) is satisfied. On the other hand, we have also 
$$
 \mathbb{E} \big( f_n^2 \,\, \mathbf{ 1}_{ \{ |f_n| \le M \} } \big) \,=\, 2 \,c  \, \big(1-\varrho_n \big) \sum_{2 \le k \le M} \, k^2 \cdot \frac{1}{\, k^2 \log k \,}\,\le \,  \frac{\,2 \,c \,M\,}{\log 2\,} \big(1-\varrho_n \big)
$$
and this shows that the condition (\ref{2.5}) is also satisfied when 
\begin{equation}
\label{4.2}
\limsup_{n \to \infty} \,\varrho_n \,=\,1\,.
\end{equation}
We conclude that, under the condition (\ref{4.2}), there is a subsequence of $f_1, f_2, \cdots \,, $ denoted again by the same symbols, for which the WLLN holds with $D_N \equiv 0$, and hereditarily:  that is, 
$$\,
\lim_{N \to \infty} \frac{1}{N} \sum_{n=1}^N f_n \,=\,0
   \quad
 \text{holds in probability}
 $$
for $f_1, f_2, \cdots \,$   and   for every one of its   subsequences. 
}
\end{example}   

\begin{remark}   
{\rm
Theorem \ref{thm2.1}  
has a   direct extension, with only very obvious notational changes, to the case where $f_1, f_2, \cdots$ take values in some  Euclidean space $\R^d$, rather than  the real line. In such an extension it does not matter whether   balls or cubes are considered in the truncation scheme   (\ref{3.3}).
}
\end{remark}



   \end{document}